\documentclass[oneside,reqno,12pt]{amsart}

\usepackage{amsmath,amssymb,amsfonts}
\usepackage{kpfonts}
\usepackage{mathtools}
\usepackage{bbm}
\usepackage{hyperref}
\usepackage[foot]{amsaddr}

\usepackage{bbm}

\newtheorem{thm}{Theorem}
\newtheorem{rem}{Remark}
\newtheorem{lemma}{Lemma}
\newtheorem{cor}{Corollary}
\newtheorem{assertion}{Proposition}
\newtheorem{example}[thm]{Example}

\newcommand{\LCM}{{\rm LCM}\,}
\newcommand{\GCD}{{\rm GCD}\,}
\newcommand{\Var}{{\rm Var}\,}

\newcommand{\E}{\mathbb{E}}

\newcommand{\me}{\mathbb{E}}
\newcommand{\mn}{\mathbb{N}}
\newcommand{\mr}{\mathbb{R}}
\newcommand{\mmp}{\mathbb{P}}

\newcommand{\tofdd}{\overset{{\rm f.d.d.}}{\longrightarrow}}
\newcommand{\tovague}{\overset{{\rm v}}{\longrightarrow}}
\newcommand{\toas}{\overset{{\rm a.s.}}{\longrightarrow}}
\newcommand{\toprobab}{\overset{{\mathbb{P}}}{\longrightarrow}}

\begin{document}

\title{Arithmetic properties of multiplicative integer-valued perturbed random walks}

\author{Victor Bohdanskyi$^1$}
\address{$^1$National Technical University of Ukraine ``Igor Sikorsky Kyiv Polytechnic Institute'', 03056 Kyiv, Ukraine. Email: vbogdanskii@ukr.net}

\author{Vladyslav Bohun$^2$}

\author{Alexander Marynych$^2$}

\author{Igor Samoilenko$^2$}
\address{$^2$Faculty of Computer Science and Cybernetics, Taras Shev\-chen\-ko National University of Kyiv, Kyiv 01601, Ukraine. Emails: vladyslavbogun@gmail.com, marynych@unicyb.kiev.ua, isamoil@i.ua}

\begin{abstract}
Let $(\xi_1, \eta_1)$, $(\xi_2, \eta_2),\ldots$ be independent identically distributed $\mathbb{N}^2$-valued random vectors with arbitrarily dependent components. The sequence $(\Theta_k)_{k\in\mn}$ defined by $\Theta_k=\Pi_{k-1}\cdot\eta_k$, where $\Pi_0=1$ and $\Pi_k=\xi_1\cdot\ldots\cdot \xi_{k}$ for $k\in\mathbb{N}$, is called a multiplicative perturbed random walk. We study arithmetic properties of the random sets $\{\Pi_1,\Pi_2,\ldots, \Pi_k\}\subset \mathbb{N}$ and $\{\Theta_1,\Theta_2,\ldots, \Theta_k\}\subset \mathbb{N}$, $k\in\mn$. In particular, we derive distributional limit theorems for their prime counts and for the least common multiple.
\end{abstract}



\maketitle

\section{Introduction}

Let $(\xi_1, \eta_1)$, $(\xi_2,\eta_2),\ldots$ be independent copies of an $\mathbb{N}^2$-valued random vector $(\xi,\eta)$ with  arbitrarily dependent components. Denote by $(\Pi_k)_{k\in\mn_0}$ (as usual, $\mn_0:=\mn\cup\{0\}$) the standard multiplicative random walk defined by
$$
\Pi_0:=1,\quad \Pi_k=\xi_1\cdot\xi_2\cdots\xi_k,\quad k\in\mn.$$
A {\it multiplicative perturbed random walk} is the sequence $(\Theta_k)_{k\in\mn}$ given by
$$
\Theta_k:=\Pi_{k-1}\cdot\eta_k,\quad k\in\mn.
$$
Note that if $\mmp\{\eta=\xi\}=1$, then $\Pi_k=\Theta_k$ for all $k\in\mn$. If $\mmp\{\xi=1\}=1$, then $(\Theta_k)_{k\in\mn}$ is just a sequence of independent copies of a random variable $\eta$. In this paper we investigate some arithmetic properties of the random sets $(\Pi_k)_{k\in\mn}$ and $(\Theta_k)_{k\in\mn}$.

To set the scene we introduce first some necessary notation. Let $\mathcal{P}$ denote the set of prime numbers. For an integer $n\in\mn$ and $p\in\mathcal{P}$, let $\lambda_p(n)$ denote the multiplicity of prime $p$ in the prime decomposition of $n$, that is,
$$
n=\prod_{p\in\mathcal{P}} p^{\lambda_p(n)}. 
$$
For every $p\in\mathcal{P}$, the function $\lambda_p:\mn\mapsto\mn_0$ is totally additive in the sense that
$$
\lambda_p(mn)=\lambda_p(m)+\lambda_p(n),\quad p\in\mathcal{P},\quad m,n\in\mathbb{N}.
$$
The set of functions $(\lambda_p)_{p\in\mathcal{P}}$ is a basic brick from which many other arithmetic functions can be constructed. For example, with
$\GCD(A)$ and $\LCM(A)$ denoting the greatest common divisor and the least common multiple of a set $A\subset\mn$, respectively, we have
$$
\GCD(A)=\prod_{p\in\mathcal{P}}p^{\min_{n\in A}\lambda_p(n)}\quad\text{and}\quad\LCM(A)=\prod_{p\in\mathcal{P}}p^{\max_{n\in A}\lambda_p(n)}.
$$

The listed arithmetic functions applied either to $A=\{\Pi_1,\ldots,\Pi_n\}$ or $A=\{\Theta_1,\ldots,\Theta_n\}$ are the main objects of investigation in the present paper. From the additivity of $\lambda_p$ we infer
\begin{equation}\label{eq:S_k(p)_def}
S_k(p):=\lambda_p(\Pi_k)=\sum_{j=1}^{k}\lambda_p(\xi_j),\quad p\in\mathcal{P},\quad k\in\mn_0,
\end{equation}
and
\begin{equation}\label{eq:T_k(p)_def}
T_k(p):=\lambda_p(\Theta_k)=\sum_{j=1}^{k-1}\lambda_p(\xi_j)+\lambda_p(\eta_k),\quad p\in\mathcal{P},\quad k\in\mn.
\end{equation}
Fix any $p\in\mathcal{P}$. Formulae~\eqref{eq:S_k(p)_def} and~\eqref{eq:T_k(p)_def}  demonstrate that $S(p):=(S_k(p))_{k\in\mn_0}$, is a standard additive random walk with the generic step $\lambda_p(\xi)$, whereas the sequence $T(p):=(T_k(p))_{k\in\mn}$, is a particular instance of an {\it additive perturbed random walk}, see~\cite{Iksanov:2016}, generated by the pair $(\lambda_p(\xi),\lambda_p(\eta))$.

\section{Main results}

\subsection{Distributional properties of the prime counts \texorpdfstring{$(\lambda_p(\xi),\lambda_p(\eta))$}{(lambda p(xi),lambda p(eta))}}

As is suggested by~\eqref{eq:S_k(p)_def} and~\eqref{eq:T_k(p)_def} the first step in the analysis of $S(p)$ and $T(p)$ should be the derivation of the joint distribution $(\lambda_p(\xi),\lambda_p(\eta))_{p\in\mathcal{P}}$. The next lemma confirms that the finite-dimensional distributions of the infinite vector $(\lambda_p(\xi),\lambda_p(\eta))_{p\in\mathcal{P}}$, are expressible via the probability mass function  of $(\xi,\eta)$. However, the obtained formulae are not easy to handle except some special cases. For $i,j\in\mathbb{N}$, put
$$
u_i:=\mmp\{\xi=i\},\quad v_j:=\mmp\{\eta=j\},\quad w_{i,j}:=\mmp\{\xi=i,\eta=j\}.
$$
\begin{lemma}\label{lem:lambda_p_joint_law}
Fix $p\in\mathcal{P}$ and nonnegative integers $(k_q)_{q\in\mathcal{P},q\leq p}$ and $(\ell_q)_{q\in\mathcal{P},q\leq p}$. Then
$$
\mmp\{\lambda_q(\xi)\geq k_q,\lambda_q(\eta)\geq \ell_q,q\in\mathcal{P},q\leq p\}=\sum_{i,j=1}^{\infty}w_{Ki,Lj},
$$
where $K:=\prod_{q\leq p,q\in\mathcal{P}}q^{k_q}$ and $L:=\prod_{q\leq p,q\in\mathcal{P}}q^{\ell_q}$.
\end{lemma}
\begin{proof}
This follows from
\begin{multline*}
\mmp\{\lambda_q(\xi)\geq k_q,\lambda_q(\eta)\geq \ell_q,q\in\mathcal{P},q\leq p\}\\
=\mmp\left\{\prod_{q\leq p,q\in\mathcal{P}}q^{k_q}\text{ divides }\xi,\prod_{q\leq p,q\in\mathcal{P}}q^{\ell_q}\text{ divides }\eta\right\}=\sum_{i,j=1}^{\infty}w_{Ki,Lj}.
\end{multline*}
Obviously, if $\xi$ and $\eta$ are independent, then
$$
\sum_{i,j=1}^{\infty}w_{Ki,Lj}=\left(\sum_{i=1}^{\infty}u_{Ki}\right)\left(\sum_{j=1}^{\infty}v_{Lj}\right).
$$
\end{proof}

We proceed with the series of examples.

\begin{example}\label{example1}
For $\alpha>1$, let $\mathbb{P}\{\xi=k\}=(\zeta(\alpha))^{-1} k^{-\alpha}$, $k\in\mathbb{N}$, where $\zeta$ is the Riemann zeta-function. Then, $(\lambda_p(\xi))_{p\in\mathcal{P}}$ are mutually independent and
$$
\mathbb{P}\{\lambda_p(\xi)\geq k\}=\sum_{i=1}^{\infty}\mathbb{P}\{\xi=p^k i\}=p^{-k\alpha},\quad k\in\mn_0,\quad p\in\mathcal{P},
$$
which means that $\lambda_p(\xi)$ has a geometric distribution on $\mn_0$ with parameter $p^{-\alpha}$.
\end{example}

\begin{example}\label{example2}
For $\beta\in(0,1)$, let $\mathbb{P}\{\xi=k\}=\beta^{k-1}(1-\beta)$, $k\in\mathbb{N}$. Then
$$
\mathbb{P}\{\lambda_p(\xi)\geq k\}=\frac{1-\beta}{\beta}\sum_{j=1}^{\infty}\beta^{p^k j}=\frac{(1-\beta)(\beta^{p^k-1})}{1-\beta^{p^k}},\quad k\in\mathbb{N}_0.
$$
\end{example}

\begin{example}\label{example3}
Let ${\rm Poi}(\lambda)$ be a random variable with the Poisson distribution with parameter $\lambda$ and put
$$
\mathbb{P}\{\xi=k\}=\mathbb{P}\{{\rm Poi}(\lambda)=k|{\rm Poi}(\lambda)\geq 1\}=(e^{\lambda}-1)^{-1}\lambda^k/k!,\quad k\in\mathbb{N}.
$$
Then
\begin{multline}
\mathbb{P}\{\lambda_p(\xi)\geq k\}=(e^{\lambda}-1)^{-1}\sum_{j=1}^{\infty}\lambda^{p^k j}/(p^k j)!\\
=\left(_0 F_{p^k}\left(;\frac{1}{p^k},\frac{2}{p^k},\ldots,\frac{p^k-1}{p^k};\left(\frac{\lambda}{p^k}\right)^{p^k}\right)-1\right),
\end{multline}
where $_0 F_{p^k}$ is the generalized hypergeometric function, see Chapter 16 in~\cite{Olver}.
\end{example}

In all examples above the distribution of $\lambda_p(\xi)$ for every fixed $p\in\mathcal{P}$, is extremely light-tailed. It is not that difficult to construct `weird' distributions where all $\lambda_p(\xi)$ have infinite expectations.

\begin{example}\label{example4}
Let $(g_p)_{p\in\mathcal{P}}$ be any probability distribution supported by $\mathcal{P}$, $g_p>0$, and $(t_k)_{k\in\mathbb{N}_0}$ any probability distribution on $\mathbb{N}$ such that $\sum_{k=1}^{\infty} kt_k=\infty$ and $t_k>0$. Define a probability distribution $\mathfrak{h}$ on $\mathcal{Q}:=\bigcup_{p\in\mathcal{P}}\{p,p^2,\ldots\}$ by
$$
\mathfrak{h}({\{p^k\}})=g_p t_k,\quad p\in\mathcal{P},\quad k\in\mathbb{N}.
$$
If $\xi$ is a random variable with distribution $\mathfrak{h}$, then
$$
\mathbb{P}\{\lambda_p(\xi)\geq k\}=g_p\sum_{j=k}^{\infty}t_j,\quad k\in\mathbb{N},\quad p\in\mathcal{P},
$$
which implies $\mathbb{E}[\lambda_p(\xi)]=g_p\sum_{k=1}^{\infty} kt_k=\infty$, $p\in\mathcal{P}$.

This example can be modified by taking $g:=\sum_{p\in\mathcal{P}}g_p<1$ and charging all points of $\mathbb{N}\setminus \mathcal{Q}$ (this set contains $1$ and all integers having at least two different prime factors) with arbitrary positive masses of the total weight $1-g$. The obtained probability distribution charges all points of $\mathbb{N}$ and still possesses the property that all $\lambda_p$'s have infinite expectations.
\end{example}

Let $X$ be a random variable taking values in $\mathbb{N}$. Since
$$
\log X=\sum_{p\in\mathcal{P}}\lambda_p(X)\log p,
$$
we conclude that $\E [(\lambda_p(X))^k]<\infty$, for all $p\in\mathcal{P}$, whenever $\E[\log^k X]<\infty$, $k\in\mathbb{N}$. It is also clear that the converse implication is false in general. When $k=1$ the inequality $\E [\lambda_p(X)]<\infty$ is equivalent to $\sum_{p\in\mathcal{P}}\E[\lambda_p(X)]\log p<\infty$. As we have seen in the above examples, checking that $\E [(\lambda_p(X))^k]<\infty$ might be a much more difficult task than proving a stronger assumption $\E[\log^k X]<\infty$. Thus, we shall mostly work under  moment conditions on $\log\xi$ and $\log\eta$.

Our standing assumption throughout the paper is
\begin{equation}\label{eq:finite_mean}
\mu_{\xi}:=\mathbb{E} [\log \xi]<\infty,
\end{equation}
which, by the above reasoning, implies $\mathbb{E} [\lambda_p(\xi)]<\infty$, $p\in\mathcal{P}$.

\subsection{Limit theorems for \texorpdfstring{$S(p)$}{S(p)} and \texorpdfstring{$T(p)$}{T(p)}}

From Donsker's invariance principle we immediately obtain the following proposition. Let $D:=D([0,\infty),\mathbb{R})$ be the Skorokhod space endowed with the standard $J_1$-topology.

\begin{assertion}\label{prop:MCLT_S}
Assume that $\E [\log^2 \xi]\in (0,\infty)$. Then,
$$
\left(\left(\frac{S_{\lfloor ut \rfloor}(p)-ut\E\lambda_p(\xi)}{\sqrt{t}}\right)_{u\geq 0}\right)_{p\in\mathcal{P}}~\Longrightarrow~((W_p(u))_{u\geq 0})_{p\in\mathcal{P}},\quad t\to\infty,
$$
on the product space $D^{\mathbb{N}}$, where, for all $n\in\mn$ and all $p_1<p_2<\cdots<p_n$, $p_i\in\mathcal{P}$, $i\leq n$, $(W_{p_1}(u),\ldots,W_{p_n}(u))_{u\geq 0}$ is an $n$-dimensional centered Wiener process with covariance matrix $C=||C_{i,\,j}||_{1\leq i,j\leq n}$ given by $C_{i,\,j}=C_{j,\,i}={\rm Cov}\,(\lambda_{p_i}(\xi),\lambda_{p_j}(\xi))$.
\end{assertion}

According to the proof of Proposition~1.3.13 in~\cite{Iksanov:2016}, see pp.~28-29 therein, the following holds true for the perturbed random walks $T(p)$, $p\in\mathcal{P}$.
\begin{assertion}\label{prop:MCLT_T}
Assume that $\E [\log^2 \xi]\in (0,\infty)$ and
\begin{equation}\label{eq:eta_p_negligible}
\lim_{t\to\infty}t^2\mathbb{P}\{\lambda_p(\eta)\geq t\}=0,\quad p\in\mathcal{P}.
\end{equation}
Then,
$$
\left(\left(\frac{T_{\lfloor ut \rfloor}(p)-ut\E\lambda_p(\xi)}{\sqrt{t}}\right)_{u\geq 0}\right)_{p\in\mathcal{P}}~\Longrightarrow~((W_p(u))_{u\geq 0})_{p\in\mathcal{P}},\quad t\to\infty,
$$
on the product space $D^{\mathbb{N}}$.
\end{assertion}
\begin{rem}
Since $\mathbb{P}\{\lambda_p(\eta)\log p\geq t\}\leq \mathbb{P}\{\log \eta \geq t\}$, the condition
\begin{equation}\label{eq:log_eta_negligible}
\lim_{t\to\infty}t^2\mathbb{P}\{\log\eta\geq t\}=0
\end{equation}
is clearly sufficient for~\eqref{eq:eta_p_negligible}.
\end{rem}

From the continuous mapping theorem under the assumptions of Proposition~\ref{prop:MCLT_T} we infer
\begin{multline}\label{eq:sup_t_p_MCLT}
\left(\left(\frac{\max_{1\leq k\leq \lfloor ut \rfloor}(T_{k}(p)-k\E\lambda_p(\xi))}{\sqrt{t}}\right)_{u\geq 0}\right)_{p\in\mathcal{P}}\\
\Longrightarrow~((\sup_{0\leq v\leq u}W_p(v))_{u\geq 0})_{p\in\mathcal{P}},\quad t\to\infty,
\end{multline}
see Proposition 1.3.13 in~\cite{Iksanov:2016}.

Formula~\eqref{eq:sup_t_p_MCLT}, for a fixed $p\in\mathcal{P}$, belongs to the realm of limit theorems for the maximum of a single additive perturbed random walk. This circle of problems is well-understood, see Section 1.3.3 in~\cite{Iksanov:2016} and~\cite{Iksanov+Pilipenko+Samoilenko:2017}, in the situation when the underlying additive standard random walk is {\it centered} and attracted to a stable L\'{e}vy process. In our setting the perturbed random walks $(T_k(p))_{k\in\mathbb{N}}$ and $(T_k(q))_{k\in\mathbb{N}}$ are dependent whenever $p,q\in\mathcal{P}$, $p\neq q$, which make derivation of the joint limit theorems harder and leads to various asymptotic regimes.

Note that~\eqref{eq:eta_p_negligible} implies $\E [\lambda_p(\eta)]<\infty$ and~\eqref{eq:log_eta_negligible} implies $\E [\log \eta]<\infty$. Theorem~\ref{main1} below tells us that under such moment conditions and assuming also $\E [\log^2 \xi]<\infty$ the maxima $\max_{1\leq k\leq n}\,T_{k}(p)$, $p\in\mathcal{P}$, of {\it noncentered} perturbed random walks $T(p)$ have the same behavior as $S_{n}(p)$, $p\in\mathcal{P}$ as $n\to\infty$.

\begin{thm}\label{main1}
Assume that $\me [\log^2 \xi]<\infty$ and $\me [\lambda_p(\eta)]<\infty$, $p\in\mathcal{P}$. Suppose further that
\begin{equation}\label{eq:full_support_assumption}
\mathbb{P}\{\xi\text{ is divisible by }p\}=\mathbb{P}\{\lambda_p(\xi)>0\}>0,\quad p\in\mathcal{P}.
\end{equation}
Then, as $t\to\infty$,
\begin{equation}\label{eq:main1_claim}
\left(\left(\frac{\max_{1\leq k\leq \lfloor tu\rfloor}\,T_{k}(p)-\me [\lambda_p(\xi)]tu}{t^{1/2}}\right)_{u\geq 0}\right)_{p\in\mathcal{P}}~\tofdd~((W_p(u))_{u\geq 0})_{p\in\mathcal{P}}.
\end{equation}
\end{thm}
\begin{rem}
If~\eqref{eq:full_support_assumption} holds only for some $\mathcal{P}_0\subseteq\mathcal{P}$, then~\eqref{eq:main1_claim} holds with $\mathcal{P}_0$ instead of $\mathcal{P}$.
\end{rem}

In the next result we shall assume that $\eta$ dominates $\xi$ in a sense that the asymptotic behavior of $\max_{1\leq k\leq n}T_{k}(p)$ is regulated by the perturbations $(\lambda_p(\eta_k))_{k\leq n}$ for all $p\in\mathcal{P}_0$, where $\mathcal{P}_0$ is a finite subset of prime numbers and those $p$'s dominate all other primes.

\begin{thm}\label{main11}
Assume~\eqref{eq:finite_mean}. Suppose further that there exists a finite set $\mathcal{P}_0\subseteq \mathcal{P}$, $d:=|\mathcal{P}_0|$, such that the distributional tail of $(\lambda_p(\eta))_{p\in\mathcal{P}_0}$ is regularly varying at infinity in the following sense. For some positive function $(a(t))_{t>0}$ and a  measure $\nu$ satisfying $\nu(\{x\in \mathbb{R}^d:\|x\|\geq r\})=c\cdot r^{-\alpha}$, $c>0$, $\alpha\in (0,1)$, it holds
\begin{equation}\label{eq:reg_var}
t\mathbb{P}\{(a(t))^{-1}(\lambda_p(\eta))_{p\in\mathcal{P}_0}\in\cdot\}~\tovague~ \nu(\cdot),\quad t\to\infty,
\end{equation}
on the space of locally finite measures on $(0,\infty]^d$ endowed with the vague topology. Finally, suppose $\mathbb{E}[\lambda_p(\eta)]<\infty$, for $p\in\mathcal{P}\setminus \mathcal{P}_0$. Then
\begin{equation}\label{eq:exteme_dominates1}
\left(\left(\frac{\max_{1\leq k\leq \lfloor tu\rfloor}\,T_{k}(p)}{a(t)}\right)_{u\geq 0}\right)_{p\in\mathcal{P}_0}~\tofdd~(M_p(u))_{u\geq 0})_{p\in\mathcal{P}_0},\quad t\to\infty,
\end{equation}
where $(M_p(u))_{u\geq 0})_{p\in\mathcal{P}_0}$ is a multivariate extreme process defined by
\begin{equation}\label{eq:extreme_def}
(M_p(u))_{p\in\mathcal{P}_0}=\sup_{k:\,t_k\leq u} y_k,\quad u\geq 0.
\end{equation}
Here the pairs $(t_k,y_k)$ are the atoms of a Poisson point process on $[0,\infty)\times (0,\infty]^d$ with the intensity measure $\mathbb{LEB}\otimes \nu$ and the supremum is taken coordinatewise. Moreover,
\begin{equation}\label{eq:exteme_dominates2}
\left(\left(\frac{\max_{1\leq k\leq \lfloor tu\rfloor}\,T_{k}(p)}{a(t)}\right)_{u\geq 0}\right)_{p\in\mathcal{P}\setminus \mathcal{P}_0}~\tofdd~0,\quad t\to\infty.
\end{equation}
\end{thm}

\subsection{Limit theorems for the \texorpdfstring{$\LCM$}{LCM}}

The results from the previous section will be applied below to the analysis of
$$
\mathcal{\Pi}_n:=\LCM(\{\Pi_1,\Pi_2,\ldots, \Pi_n\})\quad\text{and}\quad  \mathcal{\Theta}_n:=\LCM(\{\Theta_1, \Theta_2,\ldots, \Theta_n\}).
$$
A moment's reflection shows that the analysis of $\mathcal{\Pi}_n$ is trivial. Indeed, by definition, $\Pi_{n-1}$ divides $\Pi_n$ and thereupon $\mathcal{\Pi}_n=\Pi_n$ for $n\in\mn$. Thus, assuming that $\sigma_{\xi}^2:=\Var (\log \xi)\in (0,\infty)$, an application of the Donsker functional limit theorem  yields
\begin{equation}\label{eq:1}
\Big(\frac{\log \mathcal{\Pi}_{\lfloor tu\rfloor}-\mu_{\xi} tu}{ t^{1/2}}\Big)_{u\geq 0}~\Longrightarrow~(\sigma_{\xi} W(u))_{u\geq 0},\quad t\to\infty,
\end{equation}
on the Skorokhod space $D$, where $(W(u))_{u\geq 0}$ is a standard Brownian motion.

A simple structure of the sequence $(\mathcal{\Pi}_n)_{n\in\mn}$ breaks down completely upon introducing the perturbations $(\eta_k)$, which makes the analysis of $(\mathcal{\Theta}_n)$ a much harder problem. For instance, it contains as a special case the problem of studying the $\LCM$ of an independent sample, which is itself highly non-trivial. Note that
$$
\log \mathcal{\Theta}_n=\log \prod_{p\in\mathcal{P}} p^{\max_{1\leq k\leq n}\,(\lambda_p(\xi_1)+\ldots +\lambda_p(\xi_{k-1})+\lambda_p(\eta_k))}=\sum_{p\in\mathcal{P}}\max_{1\leq k\leq n} T_k(p)\log p,
$$
which shows that the asymptotic of $\mathcal{\Theta}_n$ is intimately connected with the behavior of $\max_{1\leq k\leq n}T_k(p)$, $p\in\mathcal{P}$.

As one can guess from Theorem~\ref{main1} in a `typical' situation relation~\eqref{eq:1} holds with $\log \mathcal{\Theta}_{\lfloor tu\rfloor}$ replacing $\log \mathcal{\Pi}_{\lfloor tu\rfloor}$. The following heuristics suggest the right form of assumptions ensuring that perturbations $(\eta_k)_{k\in\mn}$ have an asymptotically negligible impact on $\log \mathcal{\Theta}_n$. Take a prime $p\in\mathcal{P}$. Its contribution to $\log \mathcal{\Theta}_n$ (up to a factor $\log p$) is $\max_{1\leq k\leq n}T_k(p)$. According to Theorem~\ref{main1}, this maximum is asymptotically the same as $S_n(p)$. However, as $p$ gets large, the mean $\E [\lambda_p(\xi)]$ of the random walk $S_{n-1}(p)$ becomes small because of the identity
$$
\sum_{p\in\mathcal{P}}\E[\lambda_p(\xi)] \log p = \E [\log \xi]<\infty.
$$
Thus, for large $p\in\mathcal{P}$, the remainder $\max_{1\leq k\leq n}T_k(p)-S_{n-1}(p)$ can, in principle, become larger than $S_{n-1}(p)$ itself if the tail of $\lambda_p(\eta)$ is sufficiently heavy. In order to rule out such a possibility, we introduce the following deterministic sets:
\begin{equation}\label{eq:p_1_p_2_defs}
\mathcal{P}_1(n):=\{p\in\mathcal{P}: \mathbb{P}\{\lambda_p(\xi)>0\}\geq n^{-1/2}\}\quad\text{and}\quad \mathcal{P}_2(n):=\mathcal{P}\setminus \mathcal{P}_1(n),
\end{equation}
and bound the rate of growth of $\max_{1\leq k\leq n}\lambda_p(\eta_k)$ for all $p\in\mathcal{P}_2(n)$. It is important to note that under the assumption~\eqref{eq:full_support_assumption} it holds
$$
\lim_{n\to\infty}\min \mathcal{P}_2(n)=\infty.
$$
Therefore, if $\E [\log X]<\infty$ for some random variable $X$, then the relation
\begin{equation}\label{eq:p_2_explanation}
\lim_{n\to\infty}\sum_{p\in\mathcal{P}_2(n)} \E [\lambda_p(X)]\log p=0,
\end{equation}
holds true.
\begin{thm}\label{main2}
Assume~$\E [\log^2 \xi]<\infty$, $\me [\log \eta]<\infty$, \eqref{eq:full_support_assumption} and the following two conditions
\begin{equation}\label{eq:second_moment_diff}
\sum_{p\in\mathcal{P}}\E \left[((\lambda_p(\eta)-\lambda_p(\xi))^{+})^2\right]\log p<\infty
\end{equation}
and
\begin{equation}\label{eq:main2_eta_negligible}
\sum_{p\in\mathcal{P}_2(n)}\E [(\lambda_p(\eta)-\lambda_p(\xi))^{+}] \log p=o(n^{-1/2}),\quad n\to\infty.
\end{equation}
Then
\begin{equation}\label{eq:2}
\left(\frac{\log \mathcal{\Theta}_{\lfloor tu\rfloor}-\mu_{\xi} tu}{t^{1/2}}\right)_{u\geq 0}~\tofdd~(\sigma_{\xi} W(u))_{u\geq 0},\quad t\to\infty,
\end{equation}
where $\mu_{\xi}=\me [\log\xi]<\infty$, $\sigma_{\xi}^2=\Var [\log \xi]$ and $(W(u))_{u\geq 0}$ is a standard Brownian motion.
\end{thm}

\begin{rem}
If $\E [\log^2 \eta]<\infty$, then~\eqref{eq:second_moment_diff} holds true. Indeed, since we assume $\E [\log^2 \xi]<\infty$,
\begin{align*}
&\hspace{-1cm}\E \left[\sum_{p\in\mathcal{P}}((\lambda_p(\eta)-\lambda_p(\xi))^{+})^2\log p\right]\leq
\E \left[\sum_{p\in\mathcal{P}}(\lambda^2_p(\eta)+\lambda_p^2(\xi))\log p\right]\\
&\leq \E \left[\left(\sum_{p\in\mathcal{P}}\lambda_p(\eta)\log p\right)^2\right]+\E \left[\left(\sum_{p\in\mathcal{P}}\lambda_p(\xi)\log p\right)^2\right]\\
&=\E [\log^2 \eta]+\E [\log^2 \xi]<\infty.
\end{align*}
The condition~\eqref{eq:main2_eta_negligible} can be replaced by a stronger one which only involves distribution of $\eta$, namely
\begin{equation}\label{eq:main2_eta_negligible_alt}
\sum_{p\in\mathcal{P}_2(n)}\E [\lambda_p(\eta)] \log p=o(n^{-1/2}),\quad n\to\infty.
\end{equation}
Taking into account~\eqref{eq:p_2_explanation} and the fact that $\E [\log \eta]<\infty$, the assumption~\eqref{eq:main2_eta_negligible_alt} is nothing else but a condition of the speed of convergence of the series 
$$
\sum_{p\in\mathcal{P}} \E[\lambda_p(\eta)]\log p=\E [\log \eta].
$$
\end{rem}

\begin{example}\label{eq:explicit}
In the settings of Example~\ref{example1} let $\xi$ and $\eta$ be arbitrarily dependent with
$$
\mathbb{P}\{\xi=k\}=\frac{1}{\zeta(\alpha)k^{\alpha}},\quad \mathbb{P}\{\eta=k\}=\frac{1}{\zeta(\beta)k^{\beta}},\quad k\in\mathbb{N},
$$
for some $\alpha,\beta>1$. Note that $\E [\log^2 \xi]<\infty$ and $\E [\log^2 \eta]<\infty$. Direct calculations show that
\begin{align*}
\mathcal{P}_1(n)&=\{p\in\mathcal{P}: p^{-\alpha}\geq n^{-1/2}=\{p\in\mathcal{P}: p\leq n^{1/(2\alpha)}\},\\
\mathcal{P}_2(n)&=\{p\in\mathcal{P}: p > n^{1/(2\alpha)}\}.
\end{align*}
From the chain of relations
$$
\E [\lambda_p(\eta)]=\sum_{j\geq 1}\mathbb{P}\{\lambda_p(\eta)\geq j\}=\sum_{j\geq 1}p^{-\beta j}=\frac{p^{-\beta}}{1-p^{-\beta}}\leq 2p^{-\beta},
$$
we obtain that
\begin{multline*}
\sum_{p\in\mathcal{P}_2(n)}\E [\lambda_p(\eta)] \log p\leq 2\sum_{p\in\mathcal{P},p>n^{1/(2\alpha)}}p^{-\beta} \log p\\
\sim~2\int_{n^{1/(2\alpha)}}^{\infty}x^{-\beta}\log x \frac{{\rm d}x}{\log x}=\frac{2n^{(1-\beta)/(2\alpha)}}{\beta-1},\quad n\to\infty,
\end{multline*}
where we have used the prime number theorem for the asymptotic equivalence. Thus,~\eqref{eq:main2_eta_negligible_alt} holds if
$$
\frac{1}{2}+\frac{1-\beta}{2\alpha}<0~\Longleftrightarrow~\alpha+1<\beta.
$$
\end{example}

In the setting of Theorem~\ref{main11} the situation is much simpler in a sense that almost no extra assumptions are needed to derive a limit theorem for $\mathcal{\Theta}_n$.

\begin{thm}\label{main21}
Under the same assumptions as in Theorem~\ref{main11} and assuming additionally that
\begin{equation}\label{eq:eta_without_P0}
\sum_{p\in\mathcal{P}\setminus\mathcal{P}_0}  \E[\lambda_p(\eta)]\log p<\infty,
\end{equation}
it holds
\begin{equation}\label{eq:LCM_to_extreme}
\left(\frac{\log \mathcal{\Theta}_{\lfloor tu\rfloor}}{a(t)}\right)_{u\geq 0}~\tofdd~\left(\sum_{p\in\mathcal{P}_0} M_p(u)\log p\right)_{u\geq 0},\quad t\to\infty.
\end{equation}
\end{thm}
Note that it is allowed to take in Theorem~\ref{main21} $\xi=1$, which yields the following limit theorem for the $\LCM$ of an independent integer-valued random variables.
\begin{cor}
 Under the same assumptions on $\eta$ as in Theorem~\ref{main11} it holds
$$
\left(\frac{\log \LCM(\eta_1,\eta_2,\ldots,\eta_{\lfloor tu\rfloor})}{a(t)}\right)_{u\geq 0}~\tofdd~\left(\sum_{p\in\mathcal{P}_0} M_p(u)\log p\right)_{u\geq 0},\quad t\to\infty.
$$
\end{cor}

\begin{rem}
The results presented in Theorems~\ref{main2} and~\ref{main21} is a contribution to a popular topic in probabilistic number theory, namely, the asymptotic analysis of the $\LCM$ of various random sets. For random sets comprised of independent random variables uniformly distributed on $\{1,2,\ldots,n\}$ this problem has been addressed in~\cite{BosMarRasch:2019,BurIksMar:2022,Fernandez+Fernandez:2021,Hilberdink+Toth:2016,Kim}. Some models with a more sophisticated dependence structure have been studied~\cite{AlsKabMar:2019} and~\cite{KabMarRasch:2023}.

\end{rem}

\section{Limit theorems for coupled perturbed random walks}

Theorems~\ref{main1} and~\ref{main11} will be derived from general limit theorems for the maxima of arbitrary additive perturbed random walks indexed by some parameters ranging in a countable set in the situation when the underlying additive standard random walks are positively divergent and attracted to a Brownian motion.

Let $\mathcal{A}$ be a countable or finite set of real numbers and
$$
((X(r), Y(r)))_{r\in\mathcal{A}},\quad ((X(r),Y(r)))_{r\in \mathcal{A}},\ldots
$$
be independent copies of an $\mr^{2\times |\mathcal{A}|}$ random vector $(X(r),Y(r))_{r\in \mathcal{A}}$ with arbitrarily dependent components. 
For each $r\in\mathcal{A}$, the sequence $(S^\ast_k(r))_{k\in\mn_0}$ given by
$$
S^\ast_0(r):=0,\quad S^\ast_k(r):=X_1(r)+\ldots+X_k(r),\quad k\in\mn,
$$
is an additive standard random walk. For each $r\in\mathcal{A}$, the sequence
$(T^\ast_k(r))_{k\in\mn}$ defined by
$$
T^\ast_k(r):=S^\ast_{k-1}(r)+Y_k(r),\quad k\in\mn,
$$
is an additive perturbed random walk. The sequence $((T^\ast_k(r))_{k\in\mn})_{r\in\mathcal{A}}$ is a collection of (generally) dependent additive perturbed random walks.

\begin{assertion}\label{main3}
Assume that, for each $r\in\mathcal{A}$, $\mu(r):=\me [X(r)]\in (0,\infty)$, $\Var [X(r)]\in [0,\infty)$ and $\me [Y(r)]<\infty$. Then
\begin{equation}\label{eq:main3_FLT}
\left(\left(\frac{\max_{1\leq k\leq \lfloor tu\rfloor}\,T^\ast_{k}(r)-\mu(r)tu}{t^{1/2}}\right)_{u\geq 0}\right)_{r\in\mathcal{A}}~\tofdd~((W_r(u))_{u\geq 0})_{r\in\mathcal{A}},\quad t\to\infty,
\end{equation}
where, for all $n\in\mn$ and arbitrary $r_1<r_2<\ldots<r_n$ with $r_i\in\mathcal{A}$, $i\leq n$,  $(W_{r_1}(u),\ldots, W_{r_n}(u))_{u\geq 0}$ is an $n$-dimensional centered Wiener process with covariance matrix $C=||C_{i,\,j}||_{1\leq i,j\leq n}$ with the entries $C_{i,\,j}=C_{j,\,i}={\rm Cov}\,(X(r_i), X(r_j))$.
\end{assertion}
\begin{proof}
We shall prove an equivalent statement that, as $t\to\infty$,
$$
\left(\left(\frac{\max_{0\leq k\leq \lfloor tu\rfloor}\,T^\ast_{k+1}(r)-\mu(r)tu}{t^{1/2}}\right)_{u\geq 0}\right)_{r\in\mathcal{A}}~\tofdd~((W_r(u))_{u\geq 0})_{r\in\mathcal{A}},
$$
which differs from~\eqref{eq:main3_FLT} by a shift of the subscript $k$. By the multidimensional Donsker theorem,
\begin{equation}\label{eq:donsker}
\left(\left(\frac{S^\ast_{\lfloor tu\rfloor}(r)-\mu(r)tu}{t^{1/2}}\right)_{u\geq 0}\right)_{r\in\mathcal{A}}~\Longrightarrow~\left((W_r(u))_{u\geq 0}\right)_{r\in\mathcal{A}},\quad t\to\infty,
\end{equation}
in the product topology of $D^{\mn}$. Fix any $r\in\mathcal{A}$ and write 
\begin{multline*}
\max_{0\leq k\leq \lfloor tu\rfloor}\,T^\ast_{k+1}(r)-\mu(r)tu\\
=\max_{0\leq k\leq \lfloor tu\rfloor}\,(S^\ast_k(r)-S^\ast_{\lfloor tu\rfloor}(r)+Y_{k+1}(r))+S^\ast_{\lfloor tu\rfloor}(r)-\mu(r)tu.
\end{multline*}
In view of \eqref{eq:donsker} the proof is complete once we can show that
\begin{equation}\label{eq:proof_prw_joint1}
n^{-1/2}\left(\max_{0\leq k\leq n}\,\left(S^\ast_k(r)-S^\ast_{n}(r)+Y_{k+1}(r)\right)\right)~\overset{{\mmp}}{\to}~0,\quad n\to\infty.
\end{equation}
Let $(X_0(r), Y_0(r))$ be a copy of $(X(r), Y(r))$ which is independent of $(X_k(r), Y_k(r))_{k\in\mn}$. Since the collection $$((X_1(r), Y_1(r)),\ldots, (X_{n+1}(r), Y_{n+1}(r)))$$ has the same distribution as $$((X_{n}(r), Y_{n}(r)),\ldots, (X_0(r),Y_0(r))),$$ the variable $$\max_{0\leq k\leq n}\,(S^\ast_k(r)-S^\ast_{n}(r)+Y_{k+1}(r))$$ has the same distribution as
$$
\max\big(Y_0(r), \max_{0\leq k\leq n-1}\,(-S^\ast_k(r)+Y_{k+1}(r)-X_{k+1}(r))\big).
$$

By assumption, $\me (-S^\ast_1(r))\in (-\infty, 0)$ and $\me (Y(r)-X(r))^+<\infty$. Hence, by Theorem 1.2.1 and Remark 1.2.3 in \cite{Iksanov:2016},
$$
\lim_{k\to\infty}(-S^\ast_k(r)+Y_{k+1}(r)-X_{k+1}(r))=-\infty\quad\text{a.s.}
$$
As a consequence, the a.s.\ limit
\begin{multline*}
\lim_{n\to\infty} \max\left(Y_0(r), \max_{0\leq k\leq n-1}\,(-S^\ast_k(r)+Y_{k+1}(r)-X_{k+1}(r)\right)\\
=\max\left(Y_0(r), \max_{k\geq 0}\,(-S^\ast_k(r)+Y_{k+1}(r)-X_{k+1}(r)\right)
\end{multline*}
is a.s.\ finite. This completes the proof of~\eqref{eq:proof_prw_joint1}.
\end{proof}

\begin{proof}[Proof of Theorem~\ref{main1}]
We apply Proposition~\ref{main3} with $\mathcal{A}=\mathcal{P}$, $X(p)=\lambda_p(\xi)$ and $Y(p)=\lambda_p(\eta)$. The assumption~\eqref{eq:full_support_assumption} in conjunction with $\me [\log^2 \xi]<\infty$ imply that $\me [\lambda_p(\xi)]\in (0,\infty)$ and $\Var[\lambda_p(\xi)]\in [0,\infty)$, for all $p\in\mathcal{P}$. Similarly, $\E [\lambda_p(\eta)]<\infty$ also holds.
\end{proof}

\begin{assertion}\label{prop:eta_dominatesY}
Assume $\mathbb{E}[X(r)]<\infty$, $r\in\mathcal{A}$. Assume further that there exists a finite set $\mathcal{A}_0\subseteq \mathcal{A}$, $d:=|\mathcal{A}_0|$, such that the distributional tail of $(Y(r))_{r\in\mathcal{A}_0}$ is regularly varying at infinity in the following sense. For some positive function $(a(t))_{t>0}$ and a  measure $\nu$ satisfying $\nu(\{x\in \mathbb{R}^d:\|x\|\geq r\})=c\cdot r^{-\alpha}$, $c>0$, $\alpha\in (0,1)$, it holds
\begin{equation}\label{eq:reg_varY}
t\mathbb{P}\{(a(t))^{-1}(Y(r))_{r\in\mathcal{A}_0}\in\cdot\}~\tovague~ \nu(\cdot),\quad t\to\infty,
\end{equation}
on the space of locally finite measures on $(0,\infty]^d$ endowed with the vague topology. If $\mathbb{E}[|Y(r)|]<\infty$, for $r\in\mathcal{A}\setminus \mathcal{A}_0$, then
\begin{equation}\label{eq:exteme_dominates1Y}
\left(\left(\frac{\max_{1\leq k\leq \lfloor tu\rfloor}\,T^{\ast}_{k}(r)}{a(t)}\right)_{u\geq 0}\right)_{r\in\mathcal{A}_0}~\tofdd~(M_r(u))_{u\geq 0})_{r\in\mathcal{A}_0},\quad t\to\infty,
\end{equation}
where $(M_r(u))_{u\geq 0})_{r\in\mathcal{A}_0}$ is defined as in~\eqref{eq:extreme_def}. Moreover,
\begin{equation}\label{eq:exteme_dominates2Y}
\left(\left(\frac{\max_{1\leq k\leq \lfloor tu\rfloor}\,T^{\ast}_{k}(r)}{a(t)}\right)_{u\geq 0}\right)_{r\in\mathcal{A}\setminus \mathcal{A}_0}~\tofdd~0,\quad t\to\infty.
\end{equation}
\end{assertion}
\begin{proof}
According to Corollary 5.18 in~\cite{Resnick}
$$
\left(\left(\frac{\max_{1\leq k\leq \lfloor tu\rfloor}Y_{k}(r)}{a(t)}\right)_{u\geq 0}\right)_{r\in\mathcal{A}_0}~\Longrightarrow~\left((M_r(u))_{u\geq 0}\right)_{r\in\mathcal{A}_0},\quad t\to\infty,
$$
in the product topology of $D^{\mn}$. The function $(a(t))_{t\geq 0}$ is regularly varying at infinity with index $1/\alpha>1$. Thus, by the law of large numbers, for all $r\in\mathcal{A}$,
\begin{align}
&\left(\frac{\min_{1\leq k\leq \lfloor tu\rfloor} S^{\ast}_{k-1}(r)}{a(t)}\right)_{u\geq 0}~\tofdd~0,\quad t\to\infty,\label{eq:exteme_dominates_proof1-1}\\
&\left(\frac{\max_{1\leq k\leq \lfloor tu\rfloor} S^{\ast}_{k-1}(r)}{a(t)}\right)_{u\geq 0}~\tofdd~0,\quad t\to\infty,\label{eq:exteme_dominates_proof1-2}
\end{align}
and~\eqref{eq:exteme_dominates1Y} follows from the inequalities
\begin{multline*}
\min_{1\leq k\leq \lfloor tu\rfloor} S^{\ast}_{k-1}(r) +\max_{1\leq k\leq \lfloor tu\rfloor} Y_{k}(r)\leq \max_{1\leq k\leq \lfloor tu\rfloor} T^{\ast}_{k}(r)\\
\leq \max_{1\leq k\leq \lfloor tu\rfloor} S^{\ast}_{k-1}(r) +\max_{1\leq k\leq \lfloor tu\rfloor} Y_{k}(r).
\end{multline*}
In view of~\eqref{eq:exteme_dominates_proof1-1} and~\eqref{eq:exteme_dominates_proof1-2} , to prove~\eqref{eq:exteme_dominates2Y} it suffices to check that
$$
\left(\left(\frac{\max_{1\leq k\leq \lfloor tu\rfloor}Y_{k}(r)}{a(t)}\right)_{u\geq 0}\right)~\tofdd~0,\quad t\to\infty,
$$
for every fixed $r\in\mathcal{A}\setminus\mathcal{A}_0$. This, in turn, follows from
$$
\frac{Y_{n}(r)}{n}~\toas~0,\quad n\to\infty,\quad r\in\mathcal{A}\setminus\mathcal{A}_0,
$$
which is a consequence of the assumption $\E[|Y(r)|]<\infty$, $r\in\mathcal{A}\setminus\mathcal{A}_0$ and the Borel-Cantelli lemma.

\end{proof}
\begin{proof}[Proof of Theorem~\ref{main11}]
Follows immediately from Proposition~\ref{prop:eta_dominatesY} applied with $\mathcal{A}=\mathcal{P}$, $X(p)=\lambda_p(\xi)$ and $Y(p)=\lambda_p(\eta)$.
\end{proof}

\section{Proof of Theorem~\ref{main2}}

We aim at proving that
\begin{equation}\label{eq:main2_proof1}
\frac{\sum_{p\in\mathcal{P}}\left(\max_{1\leq k\leq n} T_k(p)-S_{n-1}(p)\right)\log p}{\sqrt{n}}~\toprobab 0,\quad n\to\infty,
\end{equation}
which together with the relation
$$
\sum_{p\in\mathcal{P}}S_{n}(p)\log p = \log \Pi_{n}= \log \mathcal{\Pi}_n,\quad n\in\mathbb{N},
$$
implies Theorem~\ref{main2} by Slutskiy's lemma and~\eqref{eq:1}.

Let $(\xi_0,\eta_0)$ be an independent copy of $(\xi,\eta)$ which is also independent of $(\xi_n,\eta_n)_{n\in\mn}$. By the same reasoning as we have used in the proof of~\eqref{eq:proof_prw_joint1} we obtain
$$
(\max_{1\leq k\leq n} T_k(p)-S_{n-1}(p))_{p\in\mathcal{P}}
\overset{d}{=}\left(\max\left(\lambda_p(\eta_0),\max_{1\leq k<n}(\lambda_p(\eta_k)-\lambda_p(\xi_k)-S_{k-1}(p))\right)\right)_{p\in\mathcal{P}}.
$$
Taking into account
$$
\sum_{p\in\mathcal{P}}\lambda_p(\eta_0)\log p=\log \eta_0,
$$
we see that~\eqref{eq:main2_proof1} is a consequence of
\begin{equation}\label{eq:main2_proof2}
\frac{\sum_{p\in\mathcal{P}}\max_{1\leq k<n}\left(\lambda_p(\eta_k)-\lambda_p(\xi_k)-S_{k-1}(p)\right)^{+}\log p}{\sqrt{n}}~\toprobab 0,\quad n\to\infty,
\end{equation}
Since, for every fixed $p\in\mathcal{P}$,
\begin{equation}\label{eq:maximum_finite}
\max_{k\geq 1}\left(\lambda_p(\eta_k)-\lambda_p(\xi_k)-S_{k-1}(p)\right)^{+}<\infty\quad \text{a.s.}
\end{equation}
by assumption~\eqref{eq:full_support_assumption}, it suffices to check that, for every fixed $\varepsilon>0$,
\begin{equation}\label{eq:main2_proof3}
\lim_{M\to\infty}\limsup_{n\to\infty}\mathbb{P}\left\{\sum_{p\in\mathcal{P},p>M}\max_{1\leq k<n}\left(\lambda_p(\eta_k)-\lambda_p(\xi_k)-S_{k-1}(p)\right)^{+}\log p>\varepsilon\sqrt{n}\right\}.
\end{equation}
In order to check~\eqref{eq:main2_proof3} we divide the sum into two disjoint parts with summations over $\mathcal{P}_1(n)$ and $\mathcal{P}_2(n)$. For the first sum, by Markov's inequality, we obtain
\begin{align*}
&\hspace{-0.4cm}\mathbb{P}\left\{\sum_{p\in\mathcal{P}_1(n),p>M}\max_{1\leq k<n}\left(\lambda_p(\eta_k)-\lambda_p(\xi_k)-S_{k-1}(p)\right)^{+}\log p>\varepsilon\sqrt{n}/2\right\}\\
&\leq \frac{2}{\varepsilon\sqrt{n}}\sum_{p\in\mathcal{P}_1(n),p>M}\E\left(\max_{1\leq k<n}\left(\lambda_p(\eta_k)-\lambda_p(\xi_k)-S_{k-1}(p)\right)^{+}\right)\log p\\
&\leq \frac{2}{\varepsilon\sqrt{n}}\sum_{p\in\mathcal{P}_1(n),p>M}\log p\sum_{k\geq 1}\E\left(\lambda_p(\eta_k)-\lambda_p(\xi_k)-S_{k-1}(p)\right)^{+}\\
&= \frac{2}{\varepsilon\sqrt{n}}\sum_{p\in\mathcal{P}_1(n),p>M}\log p\sum_{j\geq 1}\mathbb{P}\{\lambda_p(\eta)-\lambda_p(\xi)=j\}\sum_{k\geq 1}\E (j-S_{k-1}(p))^{+}\\
&\leq \frac{2}{\varepsilon\sqrt{n}}\sum_{p\in\mathcal{P}_1(n),p>M}\log p\sum_{j\geq 1}j\mathbb{P}\{\lambda_p(\eta)-\lambda_p(\xi)=j\}\sum_{k\geq 0}\mathbb{P}\{S_{k}(p)\leq j\}\\
&\leq \frac{2}{\varepsilon\sqrt{n}}\sum_{p\in\mathcal{P}_1(n),p>M}\log p\sum_{j\geq 1}j\mathbb{P}\{\lambda_p(\eta)-\lambda_p(\xi)=j\}\frac{2 j}{\E [(\lambda_p(\xi)\wedge j)]},
\end{align*}
where last estimate is a consequence of Erickson's inequality for renewal functions, see Eq.~(6.5) in~\cite{Iksanov:2016}. Further, since for $p\in\mathcal{P}_1(n)$,
$$
\E [(\lambda_p(\xi)\wedge j)]\geq \mathbb{P}\{\lambda_p(\xi)\geq 1\}=\mathbb{P}\{\lambda_p(\xi)>0\}\geq n^{-1/2},
$$
we obtain
\begin{align*}
&\hspace{-1cm}\mathbb{P}\left\{\sum_{p\in\mathcal{P}_1(n),p>M}\max_{1\leq k<n}\left(\lambda_p(\eta_k)-\lambda_p(\xi_k)-S_{k-1}(p)\right)^{+}\log p>\varepsilon\sqrt{n}/2\right\}\\
&\leq \frac{4}{\varepsilon}\sum_{p\in\mathcal{P}_1(n),p>M}\log p \E \left[ ((\lambda_p(\eta)-\lambda_p(\xi))^{+})^2\right]\\
&\leq \frac{4}{\varepsilon}\sum_{p\in\mathcal{P},p>M}\log p \E \left[((\lambda_p(\eta)-\lambda_p(\xi))^{+})^2\right].
\end{align*}
The right-hand side converges to $0$, as $M\to\infty$ by~\eqref{eq:second_moment_diff}. For the sum over $\mathcal{P}_2(n)$ the derivation is simpler. By Markov's inequality
\begin{align*}
&\hspace{-0.4cm}\mathbb{P}\left\{\sum_{p\in\mathcal{P}_2(n),p>M}\max_{1\leq k<n}\left(\lambda_p(\eta_k)-\lambda_p(\xi_k)-S_{k-1}(p)\right)^{+}\log p>\varepsilon\sqrt{n}/2\right\}\\
&\leq \frac{2}{\varepsilon\sqrt{n}}\E \left[\sum_{p\in\mathcal{P}_2(n),p>M}\max_{1\leq k<n}\left(\lambda_p(\eta_k)-\lambda_p(\xi_k)-S_{k-1}(p)\right)^{+}\log p\right]\\
&\leq \frac{2n}{\varepsilon\sqrt{n}}\E \left[\sum_{p\in\mathcal{P}_2(n),p>M}\left(\lambda_p(\eta_k)-\lambda_p(\xi_k)\right)^{+}\log p\right],
\end{align*}
and the right-hand side tends to zero as $n\to\infty$ in view of~\eqref{eq:main2_eta_negligible}. The proof is complete.

\section{Proof of Theorem~\ref{main21}}
From Theorem~\ref{main11} with the aid of the continuous mapping theorem we conclude that
$$
\left(\frac{\sum_{p\in\mathcal{P}_0}\max_{1\leq k\leq \lfloor tu\rfloor}T_k(p)\log p}{a(t)}\right)_{u\geq 0}~\tofdd~\left(\sum_{p\in\mathcal{P}_0}M_p(u)\log p\right)_{u\geq 0},
$$
as $t\to\infty$. It suffices to check
\begin{equation}\label{eq:main21_proof1}
\left(\frac{\sum_{p\in\mathcal{P}\setminus\mathcal{P}_0}\max_{1\leq k\leq \lfloor tu\rfloor}T_k(p)\log p}{a(t)}\right)_{u\geq 0}~\tofdd~0,\quad t\to\infty.
\end{equation}
Since $(a(t))$ is regularly varying at infinity,~\eqref{eq:main21_proof1} follows from
\begin{equation}\label{eq:main21_proof2}
\frac{\sum_{p\in\mathcal{P}\setminus\mathcal{P}_0}\E [\max_{1\leq k\leq n}T_k(p)]\log p}{a(n)}~\to~0,\quad n\to\infty,
\end{equation}
by Markov's inequality. To check the latter note that
\begin{align*}
&\hspace{-1cm}\sum_{p\in\mathcal{P}\setminus\mathcal{P}_0}\E [\max_{1\leq k\leq n}T_k(p)]\log p\leq \sum_{p\in\mathcal{P}\setminus\mathcal{P}_0}\E [S_{n-1}(p)+\max_{1\leq k\leq n}\lambda_p(\eta_k)]\log p\\
&\leq (n-1)\sum_{p\in\mathcal{P}\setminus\mathcal{P}_0}\E [\lambda_p(\xi)]\log p+n\sum_{p\in\mathcal{P}\setminus\mathcal{P}_0}\E [\lambda_p(\eta)]\log p\\
&\leq (n-1)\E [\log\xi]+n\sum_{p\in\mathcal{P}\setminus\mathcal{P}_0}\E [\lambda_p(\eta)]\log p=O(n),\quad n\to\infty,
\end{align*}
where we have used that $\E [\log \xi]<\infty$ and the assumption~\eqref{eq:eta_without_P0}. Using that $\alpha\in (0,1)$ and $(a(t))$ is regularly varying at infinity with index $1/\alpha$, we obtain~\eqref{eq:main21_proof2}.

\section*{Acknowledgment}
The research was supported by the National Research Foundation of Ukraine (project 2020.02/0014 `Asymptotic regimes of perturbed random walks: on the edge of modern and classical probability').

\bibliographystyle{plain}
\bibliography{LCM_PRW2023}

\end{document}